\documentclass[11pt]{article}
\usepackage{amsthm}
\usepackage{amsmath}
\usepackage{amssymb}
\usepackage{framed, fullpage}

\newtheorem{theorem}               {Theorem}
\newtheorem{lemma}{Lemma}[section]

\newtheorem{claim}[lemma]{Claim}
\newtheorem{definition} [lemma] {Definition}
\newtheorem{remark}     [lemma] {Remark}

\renewenvironment{proof}[1][\proofname]{\medskip\noindent{\bfseries #1: }}{\hfill$\blacksquare$\medskip}

\title{Deterministic vs Non-deterministic Graph Property Testing}
\author{
Lior Gishboliner \thanks{School of Mathematics, Tel-Aviv University, Tel-Aviv, Israel 69978.}
\and Asaf Shapira\thanks{School of Mathematics, Tel-Aviv University, Tel-Aviv, Israel 69978. Email: {\tt asafico@tau.ac.il}. Supported in part by ISF Grant 224/11 and a Marie-Curie CIG Grant 303320.}
}

\begin{document}
\date{}

\maketitle

\begin{abstract}
A graph property ${\cal P}$ is said to be testable if one can check if a graph is close or far from satisfying ${\cal P}$ using few random local inspections.
Property ${\cal P}$ is said to be non-deterministically testable if one can supply a ``certificate'' to the fact that a graph satisfies ${\cal P}$ so that once the certificate is given its correctness can be tested. The notion of non-deterministic testing of graph properties was recently introduced by Lov\'asz and Vesztergombi \cite{LV}, who proved that (somewhat surprisingly) a graph property is testable if and only if it is non-deterministically testable. Their proof used graph limits, and so it did not supply any explicit bounds. They thus asked if one can obtain a proof of their result which will supply such bounds. We answer their question positively by proving their result using Szemer\'edi's regularity lemma.

An interesting aspect of our proof is that it highlights the fact that the regularity lemma can be interpreted as saying that all graphs can be approximated by finitely many ``template'' graphs.
\end{abstract}

\section{Introduction}\label{sec:intro}

We consider properties of finite graph, where a property of graphs is simply a family of graphs closed under isomorphism. The main focus of our paper is the following notion of efficiently checking if a graph satisfies property ${\cal P}$ or is $\epsilon$-far from satisfying it, where a graph $G$ is said to be $\epsilon$-far from satisfying ${\cal P}$ if one should add/delete at least $\epsilon n^2$ edges to turn $G$ into a graph satisfying ${\cal P}$.

\begin{definition}\label{deftest}{\bf (Testable property)}
A graph property ${\cal P}$ is called {\em testable} if there is an algorithm
$\mathcal{T}_{{\cal P}}$, called a tester, that does the following: given $\epsilon>0$ and a graph $G$, the tester $\mathcal{T}_{{\cal P}}$ samples a set $S$ of $q_{\cal P}(\epsilon)$ vertices from $G$, checks for every $i,j \in S$ whether $(i,j)\in E(G)$ and then accepts/rejects deterministically based on the subgraph of $G$ spanned by $S$. The success probability of $\mathcal{T}_{{\cal P}}$ should be at least $\frac{2}{3}$. In other words, if the input $G$ satisfies ${\cal P}$ then $\mathcal{T}_{{\cal P}}$ accepts it with probability at least $\frac{2}{3}$, and if $G$ is $\epsilon$-far from satisfying ${\cal P}$ then $\mathcal{T}_{{\cal P}}$ rejects $G$ with probability at least $\frac{2}{3}$. The function $q_{\cal P}(\epsilon)$ is called the {\em query complexity} of $\mathcal{T}_{{\cal P}}$, and does not depend on the size of the input graph.
\end{definition}

The usual definition of a property ${\cal P}$ being testable, as introduced in \cite{GGR}, allows for the algorithm to be adaptive, but as proved in \cite{AFKS,GT} one can transform any tester into a tester of the form stated in Definition \ref{deftest} with a very minor loss in the query complexity. Therefore we do not lose generality by restricting ourselves to Definition \ref{deftest}. The important point to observe about Definition \ref{deftest} is that the algorithm makes its decision solely on the basis of the random inspections it makes into the input graph $G$. In other words, the decision of the algorithm is uniquely determined by the distribution of induced subgraphs of size $q_{\cal P}(\epsilon)$ in the input graph $G$.

Starting with \cite{GGR}, the problem of characterizing the testable graph properties received a lot of attention and by now there are several general results of this type, see \cite{AFNS,LS} and the recent surveys \cite{G,RS} for more results and references on graph property testing. A drawback of these characterizations is that they are hard to state (and use). An alternative clean characterization was recently obtained by Lov\'asz and Vesztergombi \cite{LV}. To state this characterization we need a bit of notation.

A $k$-colored graph on $n$ vertices is a coloring of the edges of $K_n$ (the complete graph on $n$ vertices) using $k$ colors. Thus, a graph can be thought of as a $2$-colored graph. A property of $k$-colored graphs is again just a family of $k$-colored graphs closed under isomorphism, and it is said to be testable\footnote{We define a $k$-colored graph to be $\epsilon$-far from satisfying a property ${\cal Q}$ of $k$-colored graphs if one should modify the colors of at least $\epsilon n^2$ edges in order to turn $G$ into a $k$-colored graph satisfying ${\cal P}$.} if it is testable in the sense of Definition \ref{deftest}. A $(k,m)$-coloring of a graph $G$ is a coloring of the edges and non-edges of $G$ with the colors $\{1,\ldots,k\}$ , so that edges are colored by $\{1,...,m\}$ and non-edges are colored by $\{m+1,...,k\}$. The following is the notion of non-deterministic testing introduced in \cite{LV}.

\begin{definition}\label{maindef}{\bf (Non-deterministically testable property)}
A graph property ${\cal P}$ is called {\em non-deterministically testable} if there are integers $k,m$ and a property ${\cal Q}$ of $k$-colored graphs so that:
\begin{enumerate}
\item A graph $G$ satisfies ${\cal P}$ if and only if there is a $(k,m)$-coloring of $G$ which satisfies ${\cal Q}$.
\item ${\cal Q}$ is testable.
\end{enumerate}
\end{definition}

We are now ready to state the characterization of the testable graph properties that was obtained by Lov\'asz and Vesztergombi \cite{LV}.

\begin{theorem}\label{theo:main}{\bf (\cite{LV})}
A property ${\cal P}$ is testable if and only if it is non-deterministically testable.
\end{theorem}

Clearly any testable property is also non-deterministically testable, thus the interesting part of the above theorem is that given the fact that a property is non-deterministically testable, one can construct a standard tester for the property. Quoting \cite{LV}, ``one could say that this theorem shows that ``P=NP'' for property testing in dense
graphs''. We refer the reader to \cite{LV} for several nice illustrations showing how to apply Theorem \ref{theo:main}.

The proof of Theorem \ref{theo:main} in \cite{LV} used the machinery of graph limits. Hence, the proof was not explicit, that is, given the fact that a property ${\cal P}$ is non-deterministically testable (in the sense of Definition \ref{maindef}), it proved the existence of a standard testing algorithm for ${\cal P}$ (in the sense of Definition \ref{deftest}) but it did not supply any upper bound for the query complexity of the new tester (i.e. the function $q_{\cal P}(\epsilon)$ in Definition \ref{deftest}). Lov\'asz and Vesztergombi \cite{LV} thus asked if Theorem \ref{theo:main} can be proved using explicit arguments that will give an effective bound. Our main result in this paper gives a positive answer to their question.

Our new proof of Theorem \ref{theo:main} uses several tools related to Szemer\'edi's regularity lemma \cite{Sz}.
In Section \ref{sec:tools} we give the necessary background for applying this lemma, state some previous results as well as some preliminary lemmas that will be used in our new proof of Theorem \ref{theo:main}. As the proofs of these technical lemmas are somewhat routine we differ them to Section \ref{sec:aux}. The proof of Theorem \ref{theo:main} appears in Section \ref{sec:proof}. As our proof applies the regularity lemma, although the bounds it supplies for $q_{\cal P}(\epsilon)$ are explicit, they are rather weak ones, given by Tower-type function of $\epsilon$. Therefore, we will not keep track of the exact dependence of $q_{\cal P}(\epsilon)$ on $\epsilon$. Finally, as we mentioned in the abstract, we believe that our proof gives a nice illustration of the fact that the regularity lemma implies that all graph can be ``approximated'' using only a finitely many template graphs. In fact, this intuition is the main idea behind the proof.

\section{Tools and Preliminary Lemmas}\label{sec:tools}

As mentioned earlier, our proof of Theorem \ref{theo:main} will apply various tools related to Szemer\'edi's regularity lemma \cite{Sz}. We will start with the basic definitions, then state some previous results that we will use (Theorem \ref{test} and Lemmas \ref{close} and \ref{regular}) and then state some technical lemmas that we will need for the proof (Lemmas \ref{aprox}, \ref{differ} and \ref{chop2}). The proofs of these technical lemmas appear in Section \ref{sec:aux}. Here, and throughout the paper, when we write $x=y \pm z$ we mean that $y-z \leq x \leq y+z$.

Given two disjoint vertex sets $U,V$ we use $E(U,V)$ to denote the set of edges connecting $U$ to $V$ and set $d(U,V)=|E(U,V)|/|U||V|$ to be the density of the bipartite graph between $U$ and $V$. The basic notion of a regular pair is the following.

\begin{definition}\label{pair1}{\bf (Regular pair)}
Suppose $U,V$ are disjoint vertex sets in a graph and let $\gamma \in (0,1)$. The pair $(U,V)$ is said to be $\gamma$-regular if for every two subsets $U' \subseteq U$, $V' \subseteq V$ satisfying
$|U'| \geq \gamma|U|$, $|V'| \geq \gamma|V|$ the inequality $\left|d(U,V)-d(U',V') \right| \leq \gamma$ holds.
\end{definition}

A $\gamma$-regular pair can/should be thought of as behaving almost like a random bipartite graph of the same density.
A partition $V_1,...,V_r$ of the vertex set of a graph is called an equipartition if
$\left| |V_i|-|V_j| \right| \leq 1$ for every $1 \leq i < j \leq r$. The order of a partition $V_1,...,V_r$ is the number of parts in it (i.e. the integer $r$).

\begin{definition}\label{partition1}{\bf (Regular equipartition)}
An equipartition $V_1,...,V_r$ of the vertices of a graph is $\gamma$-regular if all but at most $\gamma\binom{r}{2}$ of the pairs $(V_i,V_j)$ are $\gamma$-regular.
\end{definition}

We now define a graph property of having a $\gamma$-regular equipartition with a predefined set of densities.

\begin{definition} {\bf (Regularity instance)}
A graph {\em regularity instance} $R$ is given by a regularity parameter $\gamma$, an integer $r$ (the order of $R$), a set of densities ${\eta_{i,j}}$ where $1\leq i<j\leq r $, and a set of non-regular pairs $\bar{R}$ of size at most $\gamma\binom{r}{2}$. A graph $G$ is said to satisfy $R$ if G has an equipartition $V_1,...,V_r$ such that for every $(i,j) \notin \bar{R}$ the pair $(V_i,V_j)$ is $\gamma$-regular and satisfies $d (V_i,V_j)=\eta_{i,j}$.
\end{definition}

A key element in the proof of Theorem \ref{theo:main} is the following result, which follows from the main results of \cite{AFNS,FN}. It allows us to test how close a graph is to satisfying a given regularity instance.

\begin{theorem}\label{test}{\bf (\cite{AFNS,FN})}
Let $R$ be a graph regularity instance, let $\epsilon_2 > \epsilon_1 > 0$ and let $p < \frac{1}{2}$. Then there is a tester $\mathcal{T}=\mathcal{T}_{\ref{test}}(R,\epsilon_1,\epsilon_2,p)$ that distinguishes graphs that are $\epsilon_1$-close to satisfying $R$ from graphs that are $\epsilon_2$-far from satisfying $R$, with success probability at least $1-p$.

Furthermore, the query-complexity of $T$ depends only on $R$, $\epsilon_1$, $\epsilon_2$ and $p$ (and not on the input graph) and can be expressed as an explicit function of these parameters.
\end{theorem}

We note that the arguments used in \cite{AFNS,FN} to prove the above result all relied heavily on the regularity lemma. Therefore, the bounds they give have a very poor (yet explicit) Tower-type dependence on the input parameters.

The second result we will need is Corollary 3.8 from \cite{AFNS}.

\begin{lemma}\label{close}{\bf (\cite{AFNS})}
Let $R$ be a graph regularity instance of order $r$, regularity parameter $\gamma$, densities ${\eta_{i,j}}$ and a set of non-regular pairs $\bar{R}$ . Suppose that a graph $G$ has an equipartition ${V_1,...,V_r}$ such that for every $(i,j) \notin \bar{R}$ the pair $(V_i,V_j)$ is $\gamma$-regular\footnote{Actually, Corollary 3.8 in \cite{AFNS} only needs to assume that $(V_i,V_j)$ is
$\left(\gamma + \frac{\gamma^{2}\epsilon}{50}\right)$-regular.} and satisfies $d(V_i,V_j)=\eta_{i,j} \pm \frac{\gamma^{2}\epsilon}{50}$.
Then $G$ is $\epsilon$-close to satisfying $R$.
\end{lemma}

We now turn to consider $k$-colored graphs. We first generalize the definitions of a regular pair, regular equipartition and regularity instance, to the more general setting of $k$-colored graphs. We start with the following notation: Suppose $U,V$ are two disjoint vertex sets in a $k$-colored graph. We use\footnote{Here, and throughout the paper, we always use $\ell$ as a superscript and never as an exponent. So $x^{\ell}$ should read ``$x$ superscript $\ell$'' not ``$x$ to the power $\ell$''.} $d^{\ell}(U,V)$ to denote the density of edges of color $\ell$ between $U$ and $V$, that is $d^{\ell}(U,V)|=E^{\ell}(U,V)|/|A||B|$, where $E^{\ell}(U,V)$ is the set of edges with color $\ell$ that connect $U$ to $V$. In case there is more than one graph, we use $d^{\ell}_G(U,V)$ to denote the density of edges colored by $\ell$ between $U,V$ in the $k$-colored graph $G$.

\begin{definition}\label{pair2} {\bf (Regular pair in a $k$-colored graph)}
Suppose $U,V$ are disjoint vertex sets in a $k$-colored graph. The pair $(U,V)$ is $\gamma$-regular if for every two subsets $U' \subseteq U$, $V' \subseteq V$ satisfying
$|U'| \geq \gamma|U|$, $|V'| \geq \gamma|V|$, and for every $1 \leq \ell \leq k$, the inequality $\left|d^{\ell}(U,V)-d^{\ell}(U',V') \right| \leq \gamma$ holds.
\end{definition}

\begin{definition}\label{partition2} {\bf (Regular equipartition in a $k$-colored graph)}
An equipartition $V_1,...,V_r$ of the vertices of a $k$-colored graph is $\gamma$-regular if all but at most $\gamma\binom{r}{2}$ of the pairs $(V_i,V_j)$ are $\gamma$-regular.
\end{definition}

\begin{definition} {\bf ($k$-colored regularity instance)}
A $k$-colored regularity instance $R$ is given by a regularity parameter $\gamma$, an integer $r$ (the order of $R$), a set of densities ${\eta_{i,j}^{\ell}}$ where $1\leq i<j\leq r $ and $1\leq \ell \leq k $, and a set of non-regular pairs $\bar{R}$ of size at most $\gamma\binom{r}{2}$. A $k$-colored graph $G$ is said to satisfy $R$ if G has an equipartition $V_1,...,V_r$ such that for every $(i,j) \notin \bar{R}$ the pair $(V_i,V_j)$ is $\gamma$-regular and satisfies $d^{\ell} (V_i,V_j)=\eta_{i,j}^{\ell}$ for every $1 \leq \ell \leq k$.
\end{definition}

The Regularity Lemma for $k$-colored graphs states that every $k$-colored graph has a $\gamma$-regular equipartition whose order can be bounded by a function of $\gamma$ and $k$. It can be formulated in terms of regularity instances in the following way.

\begin{lemma}\label{regular} {\bf (Regularity lemma for $k$-colored graphs)}
For every $\gamma >0$ and integers $t$ and $k$, there exists $T=T_{\ref{regular}}(\gamma,t,k)$ so that every $k$-colored graph satisfies some $k$-colored regularity instance of order at least $t$ and at most $T$, and regularity parameter $\gamma$.
\end{lemma}

Note that the usual regularity lemma is the special case of the $k$-colored regularity lemma with $k=2$. The proof of the $k$-colored version requires a minor adaptation of the proof of the standard regularity lemma. See \cite{KSim} for the details.

Having described the known results that will be used in the proof of Theorem \ref{theo:main}, we now turn to state the additional technical lemmas we shall rely on. We start with a lemma that allows one to approximate the number of copies of small $k$-colored graphs in a $k$-colored graph which satisfies a given regularity instance.

\begin{definition}\label{icdef1}{\bf ($IC(B,W,\sigma)$)}
Suppose $B$ is a $k$-colored graph on vertices $[q]$ in which $(i,j)$ is colored by $c(i,j)$. Suppose $W=\{{\eta_{i,j}^{\ell} : 1 \leq i<j \leq q,~ 1 \leq \ell \leq k}\}$ are densities and $\sigma:[q]\rightarrow[q]$ is a permutation. Define:
\begin{equation*}
IC(B,W,\sigma)=\prod_{i<j}{\eta^{c(\sigma(i),\sigma(j))}_{i,j}}
\end{equation*}
\end{definition}

\begin{definition}\label{icdef2}{\bf ($IC(B,W)$)}
Suppose $B$ and $W$ are as in the previous definition. Define:
\begin{equation*}
IC(B,W)=\frac{1}{Aut(B)} \sum_{\sigma}{IC(B,W,\sigma)}
\end{equation*}
where $Aut(B)$ is the number of automorphisms of $B$, that is, the number of injections $\phi: V(B) \mapsto V(B)$ that preserve the color of the edges.
\end{definition}

\begin{definition}\label{icdef3}{\bf ($IC(B,R)$)}
Let $R$ be a $k$-colored regularity instance of order $r$ and densities $\{{\eta_{i,j}^{\ell} : 1 \leq i<j \leq r,~ 1 \leq \ell \leq k}\}$. Let $B$ be a $k$-colored graph on the vertex set $[q]$. For every $A \subseteq [r]$ of size $q$ put $W(A)=\{ \eta^{\ell}_{i,j} : i,j\in A, ~1 \leq \ell \leq k \}$. Define:
\begin{equation*}
IC(B,R)=\frac{1}{\binom{r}{q}}\sum_{A \subseteq [r], |A|=q}{IC(B,W(A))}
\end{equation*}
\end{definition}

\medskip

\begin{remark} It is easy to see that $IC(B,W,\sigma)$, $IC(B,W)$ and $IC(B,R)$ are quantities in $[0,1]$.
\end{remark}

To understand Definition \ref{icdef1}, consider a random $k$-colored graph whose vertices are $V_1 \cup...\cup V_q$. Suppose that the probability that the color of $(v_i,v_j)$ is $\ell$ is $\eta^{\ell}_{i,j}$ ($v_i \in V_i, v_j \in V_j$). Suppose also that $|V_1|=...|V_q|=n$. Let $B$ be a fixed $k$-colored graph on the vertices $[q]$ and let $\sigma$ be a permutation of $[q]$. What is the expected number of $q$-tuples $v_1\in V_1,...,v_q\in V_q$ which span a copy of $B$ such that $v_i$ plays the role of $\sigma(i)$? It is easy to see that this number is $IC(B,W,\sigma)n^q$ where we set
$W=\{{\eta_{i,j}^{\ell} : 1 \leq i<j \leq q,~ 1 \leq \ell \leq k}\}$. We show (in Lemma \ref{aprox1}) that for every $\delta$, if all pairs $(V_i,V_j)$ are $\gamma$-regular for some small enough $\gamma$, then the number of such $q$-tuples $v_1,...,v_q$ is
$(IC(B,W,\sigma) \pm \delta)n^q$. This fact demonstrates the almost random behavior of regular partitions.
The expression $IC(B,W)$ (in Definition \ref{icdef2}) is used to approximate the total number of $q$-tuples $v_1\in V_1,...,v_q\in V_q$ which span a copy of $B$. The expression $IC(B,R)$ (in Definition \ref{icdef3}) is used to approximate the number of copies of $B$ in a graph that satisfies the regularity instance $R$. The most general result of this sort is the following lemma.

\begin{lemma} \label{aprox}
For any $\delta>0$ and integers $k$ and $q$ there are $\gamma=\gamma_{\ref{aprox}}(\delta,q,k)$ and $t=t_{\ref{aprox}}(\delta,q,k)$ with the following property: For any $k$-colored regularity instance $R$ of order at least $t$ and regularity parameter at most $\gamma$, and for any family $\mathcal{B}$ of $k$-colored graphs on $q$ vertices, the number of copies of $k$-colored graphs $B \in \mathcal{B}$ in any $k$-colored graph on $n$ vertices satisfying $R$ is
\begin{equation*}
\left(\sum_{B \in \mathcal{B}}{IC(B,R)} \pm \delta \right) \binom{n}{q}
\end{equation*}
\end{lemma}

\bigskip

The proof of Lemma \ref{aprox} appears in Subsection \ref{subsec:aprox}. The second lemma we will need is the following.

\begin{lemma}\label{differ}
For every $\delta$ and integers $q$ and $k$ there is
$\lambda=\lambda_{\ref{differ}}(\delta,q,k)$ such that the following holds:
Let $R$ and $M$ be $k$-colored regularity instances of order $r$, and densities
\begin{equation*}
\{ \eta^{\ell}_{i,j} : 1 \leq i<j \leq r,~ 1 \leq \ell \leq k \}
\end{equation*}
and
\begin{equation*}
 \{ \mu^{\ell}_{i,j} : 1 \leq i<j \leq r,~ 1 \leq \ell \leq k \}
\end{equation*}
respectively. Let $\mathcal{B}$ be a family of $k$-colored graphs of order $q$. Suppose that $\mu^{\ell}_{i,j}=\eta^{\ell}_{i,j} \pm \lambda$. Then
\begin{equation}\label{closeinstances}
 \left|\sum_{B \in \mathcal{B}}{IC(B,R)}- \sum_{B \in \mathcal{B}}{IC(B,M)}\right| \leq \delta\;.
\end{equation}
\end{lemma}

\bigskip

The proof of Lemma \ref{differ} appears in Subsection \ref{subsec:differ}. The last ingredient we will need is the following lemma whose proof appears in Subsection \ref{subsec:chop}.

\begin{definition}\label{chopdef}{\bf (Chopping)}
Let $R$ be a graph regularity instance of order $r$, regularity parameter $\gamma$, densities $\mu_{i,j} $ and a set $\bar{R}$ of non-regular pairs. A \textbf{$(k,m)$-chopping} of $R$ is any $k$-colored regularity instance $R'$ of order $r$, regularity parameter $2\gamma$ , non-regular set $\bar{R'}=\bar{R}$ and densities $\eta^{\ell}_{i,j}$ that satisfy
\begin{equation*}
\sum_{\ell=1}^{m}{\eta^{\ell}_{i,j}}=\mu_{i,j}~~~\mbox{and}~~~\sum_{\ell=m+1}^{k}{\eta^{\ell}_{i,j}}=1-\mu_{i,j}
\end{equation*}
\end{definition}

\begin{lemma}\label{chop2}
For every $\gamma>0$ and integers $t$ and $k$, there is $n_{\ref{chop2}}(\gamma,t,k)$ such that the following holds: Suppose $R$ is a graph regularity instance of order at most $t$ and regularity parameter $\gamma$, that $R'$ is a $(k,m)$-chopping of $R$ and that $G$ is a graph satisfying $R$ with at least $n_{\ref{chop2}}(\gamma,t,k)$ vertices. Then $G$ has a $(k,m)$-coloring that satisfies $R'$.
\end{lemma}

\section{The New Proof of Theorem \ref{theo:main}}\label{sec:proof}

Consider any $\epsilon>0$. Let property ${\cal Q}$ and integers $k$ and $m$ be those from Definition \ref{maindef}, that is, so that ${\cal Q}$ is a property of $k$-colored graphs and so that a graph satisfies ${\cal P}$ if and only if it has a $(k,m)$-coloring satisfying
${\cal Q}$. Suppose ${\cal Q}$ can be tested by a tester $\mathcal{T}_{Q}$ as in Definition \ref{deftest}. Let $q=q_{\cal Q}(\frac{\epsilon}{2})$ be the query-complexity of $\mathcal{T}_{{\cal Q}}$, i.e. the number of vertices that $\mathcal{T}_{{\cal Q}}$ samples when testing if a $k$-colored graph satisfies ${\cal Q}$ or is $\frac{\epsilon}{2}$-far from satisfying it. Let $\mathcal{B}$ be the set of all $k$-colored graphs $B$ on $q$ vertices, such that when $\mathcal{T}_{{\cal Q}}$ samples a $k$-colored graph isomorphic to $B$, it accepts the input. Put
$$
t=t_{\ref{aprox}}(1/12,q,k), ~~\gamma=\gamma_{\ref{aprox}}(1/12,q,k),~~ T=T_{\ref{regular}}(\gamma/2k,t,k)
$$
and
$$
\eta=\min\left\{ \lambda_{\ref{differ}}(1/12,q,k),\frac{\epsilon \left(\frac{\gamma}{2} \right)^2}{200m} \right\}\;.
$$
Let $I$ be the set of all $k$-colored regularity instances of order at least $t$ and at most $T$, regularity parameter $\gamma$ and densities from the set $\{0,\eta,2\eta,3\eta,...,1\}$.
Observe that all the above constants, as well as $|I|$, depend only on $\epsilon$, $k$ and the properties ${\cal P}$ and ${\cal Q}$. We now arrive at the critical definition of the proof:

\begin{definition}\label{good} {\bf (Good regularity instance)}
A graph regularity instance $R$ with regularity measure $\gamma/2$ is considered \textbf{good} if it has a $(k,m)$-chopping $R'$ that satisfies:
\begin{enumerate}
\item $R' \in I $.
\item $\sum_{B \in \mathcal{B}}{IC(B,R')} \geq \frac{1}{2}$.
\end{enumerate}
We say that $R'$ is a $\textbf{witness}$ to the fact that $R$ is good. We set $GOOD$ to be the family of good regularity instances.
\end{definition}

Suppose first that the input graph has less than $n_{\ref{chop2}}(\frac{\gamma}{2},T,k)$ vertices. In this case the algorithm can just ask about all edges of $G$ and then check if $G$ satisfies property ${\cal P}$. Since $\gamma$, $T$ and $k$ are all functions of $\epsilon$, ${\cal P}$ and ${\cal Q}$, we get that so is the query complexity in this case. Hence from this point on we will assume that $n\geq n_{\ref{chop2}}(\frac{\gamma}{2},T,k)$.

The following are the key observations we will need for the proof.

\begin{claim}\label{clm:dir1}
If $G$ satisfies ${\cal P}$, then $G$ is $\frac{\epsilon}{4}$-close to satisfying some $R \in GOOD$.
\end{claim}

\begin{claim}\label{clm:dir2}
If $G$ is $\epsilon$-far from satisfying ${\cal P}$, then $G$ is $\frac{\epsilon}{2}$-far from satisfying any $R \in GOOD$.
\end{claim}

Let us first complete the proof based on these claims. We describe a randomized algorithm that distinguishes between graphs satisfying $\mathcal{P}$ and graphs that are $\epsilon$-far from satisfying $\mathcal{P}$, with success probability at least $\frac{2}{3}$, and by making a number of queries that can be bounded by a function of $\epsilon$. Put $p=\frac{1}{3|GOOD|}$.
Let $G$ be a graph on at least $n_{\ref{chop2}}(\frac{\gamma}{2},T,k)$ vertices. In order to test $G$ for property ${\cal P}$ we do the following: For every $R \in GOOD$ use $\mathcal{T}_{13}(R, \frac{\epsilon}{4},\frac{\epsilon}{2},p)$ (recall Theorem \ref{test}) to test whether $G$ is $\frac{\epsilon}{4}$-close to satisfying $R$ or $\frac{\epsilon}{2}$-far from satisfying it. If one of these tests accepts, then accept the input $G$, otherwise reject it. If $G$ satisfies ${\cal P}$ then by Claim \ref{clm:dir1} it is $\frac{\epsilon}{4}$-close to some $R \in GOOD$, and our tester accepts it with probability at least $1-p \geq \frac{2}{3}$. If $G$ is $\epsilon$-far from satisfying ${\cal P}$ then by Claim \ref{clm:dir2} it is $\frac{\epsilon}{2}$-far from satisfying any $R \in GOOD$. Our tester accepts $G$ with probability at most $|GOOD| \cdot p=\frac{1}{3}$ and so it rejects with probability at least $\frac{2}{3}$, as required. Finally, since all the parameters involved are given by explicit functions of $\epsilon$
and the properties ${\cal P}$ and ${\cal Q}$, we get via Theorem \ref{test} that the number of queries made by the tester can be bounded by an explicit function of $\epsilon$.

We now complete the proof of Theorem \ref{theo:main} by proving Claims \ref{clm:dir1} and \ref{clm:dir2}.

\begin{proof}[Proof (of Claim \ref{clm:dir1})]
Suppose $G$ satisfies ${\cal P}$. Then there exists some $(k,m)$-coloring of $G$ that satisfies ${\cal Q}$. Denote this $k$-colored graph by $H$. By Lemma \ref{regular}, $H$ satisfies some $k$-colored regularity instance $R_{1}'$ of order $t \leq r \leq T $, regularity parameter $\frac{\gamma}{2k}$ and densities
\begin{equation*}
\{ \eta^{\ell}_{i,j} : 1 \leq i<j \leq r,~ 1 \leq \ell \leq k \}.
\end{equation*}
Since $H$ satisfies ${\cal Q}$, we infer that $\mathcal{T}_{{\cal Q}}$ must accept $H$ with probability at least $\frac{2}{3}$. This means that when sampling $q$ vertices from $H$, the probability to get a $k$-colored graph isomorphic to one of the elements of $\mathcal{B}$ is at least $\frac{2}{3}$. By the choice of $\gamma$ and $t$ via Lemma \ref{aprox} we get that this probability is $\sum_{B \in \mathcal{B}}{IC(B,R_{1}')} \pm \frac{1}{12}$. Therefore
\begin{equation}\label{eq3}
\sum_{B \in \mathcal{B}}{IC(B,R_{1}')} \geq \frac{7}{12}\;.
\end{equation}

Let $V_1,...,V_r$ be an equipartition of $H$ which corresponds to $R'_1$.
We claim that $V_1,...,V_r$ is also a $\frac{\gamma}{2}$-regular equipartition of $G$. To see this let $(i,j) \notin \bar{R_{1}'}$. For every $x \in V_i, y \in V_j$, the edge $(x,y)$ is in $G$ if and only if $(x,y)$ is colored in $H$ by a color $\ell \in \{1,...,m\}$. Therefore $d_{G}(V_i,V_j)=\sum_{\ell=1}^{m}{d_{H}^{\ell}(V_i,V_j)}$. Let $U' \subseteq V_i, V' \subseteq V_j$ such that $|U'|\geq \frac{\gamma}{2}|V_i|, |V'| \geq \frac{\gamma}{2}|V_j|$. As we assume that $(V_i,V_j)$ is a $\frac{\gamma}{2k}$-regular pair in $H$, we have $|d_{H}^{\ell}(U',V')-d_{H}^{\ell}(V_i,V_j)| \leq \frac{\gamma}{2k}$ for every $1\leq \ell \leq k$. By the triangle inequality, we have
$$
|d_{G}(U',V')-d_{G}(V_i,V_j)| \leq \sum_{\ell=1}^{m}{|d_{H}^{\ell}(U',V')-d_{H}^{\ell}(V_i,V_j)|} \leq \frac{m\gamma}{2k} \leq \frac{\gamma}{2}\;.
$$
We thus infer that $G$ satisfies a regularity instance $R_1$ with order $r$, regularity parameter $\frac{\gamma}{2}$, a set of irregular pairs $\bar{R_1'}$ and densities
$\left\{\eta_{i,j} : 1 \leq i<j \leq r \right\}$ where $\eta_{i,j}=\sum_{l=1}^{m}{\eta^{\ell}_{i,j}}$.

Let $R_{2}'$ be the $k$-colored regularity instance that is obtained from $R_{1}'$ by replacing each of the densities $\eta^{\ell}_{i,j}$ with the closest integer multiples of $\eta$. Observe that we thus change each density by at most $\eta$. As we chose
$\eta \leq \lambda_{\ref{differ}}(\frac{1}{12},q,k)$, we get from Lemma \ref{differ} and (\ref{eq3}) that
\begin{equation}\label{eq4}
\sum_{B \in \mathcal{B}}{IC(B,R_{2}')} \geq \sum_{B \in \mathcal{B}}{IC(B,R_{1}')}-\frac{1}{12} \geq \frac{1}{2}\;.
\end{equation}
Denote the densities of $R_2'$ by $\mu^{\ell}_{i,j}$. Let $R_2$ be the graph regularity instance of order $r$, regularity parameter $\frac{\gamma}{2}$, densities $\{\mu_{i,j} : 1 \leq i<j \leq r \}$, where $\mu_{i,j}=\sum_{\ell=1}^{m}{\mu^{\ell}_{i,j}}$, and a set of irregular pairs $\bar{R_2'}$. By Definition \ref{chopdef} $R_2'$ is a $(k,m)$-chopping of $R_2$. Furthermore $R_2' \in I$ and we get from (\ref{eq4}) that $\sum_{B \in \mathcal{B}}{IC(B,R_{2}')} \geq \frac{1}{2}$. By Definition \ref{good} $R_2'$ is a witness to the fact that $R_2$ is good. Finally, recalling that
$\eta \leq \frac{\epsilon \left(\frac{\gamma}{2}\right)^2}{200m}$, we get that for every $i<j$ we have
\begin{equation*}
\left|\eta_{i,j}-\mu_{i,j}\right|=\left|\sum_{l=1}^{m}{\eta^{\ell}_{i,j}}-\sum_{l=1}^{m}{\mu^{\ell}_{i,j}}\right| \leq \sum_{l=1}^{m}\left|{\eta^{\ell}_{i,j}}-{\mu^{\ell}_{i,j}}\right| \leq m\eta \leq \frac{\epsilon \left(\frac{\gamma}{2}\right)^2}{200}\;.
\end{equation*}
In other words, the densities of $R_1$ and $R_2$ differ by at most
$\frac{\epsilon \left(\frac{\gamma}{2}\right)^2}{200}$. We now get via Lemma \ref{close} that $G$ is $\frac{\epsilon}{4}$-close to satisfying $R_2$.
\end{proof}

\begin{proof}[Proof (of Claim \ref{clm:dir2})]
We will prove that if $G$ is $\frac{\epsilon}{2}$-close to satisfying some $R \in GOOD$ then $G$ is $\epsilon$-close to satisfying ${\cal P}$. Suppose that $G$ is $\frac{\epsilon}{2}$-close to a graph $G'$ that satisfies some $R \in GOOD$. By Definition \ref{good} $R$ has a $(k,m)$-chopping $R'$ such that $\sum_{B \in \mathcal{B}}{IC(B,R')} \geq \frac{1}{2}$. By Lemma \ref{chop2} $G'$ has a $(k,m)$-coloring satisfying $R'$. Call this $k$-colored graph $H'$. By Lemma \ref{aprox} the probability to get a $k$-colored graph isomorphic to an element of $\mathcal{B}$ when sampling $q$ vertices from $H'$ is $\sum_{B \in \mathcal{B}}{IC(B,R')} \pm \frac{1}{12}$. Therefore this probability is at least $\frac{5}{12}$. If $H'$ was $\frac{\epsilon}{2}$-far from satisfying ${\cal Q}$ this probability would have to be at most $\frac{1}{3}$, because $\mathcal{T}_{Q}$ would have to reject $H'$ with probability at least $\frac{2}{3}$. So we infer that $H'$ is $\frac{\epsilon}{2}$-close to satisfying ${\cal Q}$. This means that $H'$ can be turned into a $k$-colored graph $H''$ that satisfies ${\cal Q}$ by changing the colors of at most $\frac{\epsilon}{2}n^2$ edges.

Construct a graph $G''$ by doing the following: For every $x,y \in V(H'')$, put an edge between $x$ and $y$ if $(x,y)$ is colored by a color $\ell \in \{1,...,m\}$ in $H''$. First, $G''$ satisfies ${\cal P}$ because $H''$ is a $(k,m)$-coloring of $G''$ which satisfies ${\cal Q}$. Furthermore, we claim that $G''$ is $\frac{\epsilon}{2}$-close to $G'$. Indeed, observe that the number of edge modifications we performed is exactly the number of pairs $(x,y)$ so that
in one of the graphs $H',H''$ the color of $(x,y)$ belonged to the set $\{1,\ldots,m\}$ while in the other it belonged to $\{m+1,\ldots,k\}$. This number is clearly bounded from above by the number of modifications made when turning $H'$ to $H''$. Since $H''$ and $H'$ differ in at most $\frac{\epsilon}{2}n^2$ edges the same thus holds for $G''$ and $G'$, implying that $G'$ is $\frac{\epsilon}{2}$-close to $G''$. Since $G$ is assumed to be $\frac{\epsilon}{2}$-close to $G'$, we infer that $G$ is $\epsilon$-close to $G''$. Since $G''$ satisfies ${\cal P}$ the proof is complete.
\end{proof}

\section{Proofs of Auxiliary Lemmas}\label{sec:aux}

\subsection{Proof of Lemma \ref{aprox}}\label{subsec:aprox}

We will need the following folklore result stating the a $q$-tuple of vertex sets that are pairwise regular have the correct number of
copies of $K_q$ (the complete graph on $q$ vertices). A detailed proof can be found in \cite{F}.

\begin{lemma}\label{aprox1graphs}
For every $\delta>0$ and $q$ there exists $\gamma'=\gamma'_{\ref{aprox1graphs}}(\delta,q)$ such that the following holds: Suppose $V_1,...,V_q$ are disjoint vertex sets in a graph,
$|V_1|=...=|V_q|=n$, and all pairs $(V_i,V_j)$ are $\gamma'$-regular. Put $IC(K_q;V_1,...,V_q)=\prod\limits_{i<j}{d(V_i,V_j)}$.
Then the number of $q$-tuples
$v_1 \in V_1,...,v_q \in V_q$ that span a copy of $K_q$ is
\begin{equation*}
(IC(K_q; V_1,...,V_q)\pm \delta)n^q
\end{equation*}
\end{lemma}

As a first step towards proving Lemma \ref{aprox} we prove a variant of Lemma \ref{aprox1graphs} for $k$-colored graphs with respect to $IC(B,W,\sigma)$. We will then obtain similar lemmas with respect to $IC(B,W)$ and $IC(B,R)$ (recall Definitions \ref{icdef1}, \ref{icdef2} and \ref{icdef3}) and then derive from them the proof of Lemma \ref{aprox}.

\begin{lemma}\label{aprox1}
For every $\delta>0$ and $q$ there exists $\gamma=\gamma_{\ref{aprox1}}(\delta,q)$ such that the following holds: Suppose $V_1,...,V_q$ are disjoint vertex sets in a $k$-colored graph, $|V_1|=...=|V_q|=n$, and all pairs $(V_i,V_j)$ are $\gamma$-regular. Put $W=\{d^{\ell}(V_i,V_j) : 1 \leq i<j \leq q , ~1 \leq \ell \leq k\}$. Then for every $k$-colored graph $B$ on the vertices $[q]$, and for any permutation $\sigma: [q] \rightarrow [q]$, the number of $q$-tuples $v_1 \in V_1,...,v_q \in V_q$ which span copy a of $B$ with $v_i$ playing the role of $\sigma(i)$ is
\begin{equation*}
(IC(B,W,\sigma) \pm \delta)n^q
\end{equation*}
\end{lemma}

\begin{proof} While the proof of Lemma \ref{aprox1graphs} can be adapted to the more general setting of Lemma \ref{aprox1} it will be easier
to reduce Lemma \ref{aprox1} to Lemma \ref{aprox1graphs}. Set $\gamma=\gamma_{\ref{aprox1}}(\delta,q)=\gamma'_{\ref{aprox1graphs}}(\delta,q)$ and suppose $(V_i,V_j)$ is $\gamma$-regular for every $1\leq i<j \leq q$. We call a $q$-tuple $v_1\in V_1,...,v_q\in V_q$ {\bf proper}, if $v_1,...,v_q$ span a copy of $B$ with $v_i$ playing the role of $\sigma(i)$.

We denote by $c(i,j)$ the color of the edge $(i,j)$ in $B$. Let $E_{i,j}$ be the set of edges connecting a vertex in $V_i$ to a vertex in $V_j$ whose color is $c(\sigma(i),\sigma(j))$. If $v_1 \in V_1,...,v_q \in V_q$ is proper, then the color of $(v_i,v_j)$ is $c(\sigma(i), \sigma(j))$. We see that the edges in $E_{i,j}$ are the only edges between $V_i$ and $V_j$ that can "participate" in a proper $q$-tuple. Define a $q$-partite graph $S$ with vertex sets $V_1,...,V_q$, in which the edges between $V_i$ and $V_j$ are $E_{i,j}$. A $q$-tuple $v_1\in V_1,...,v_q\in V_q$ is proper if and only if it spans a copy of $K_q$ in $S$. So in order to prove Lemma \ref{aprox1} it is enough to show that the number of copies of $K_q$ in $S$ is $(IC(B,W,\sigma) \pm \delta)n^q$. By Lemma \ref{aprox1graphs}, the number of copies of $K_q$ in $S$ is
$(IC(K_q;V_1,...,V_q) \pm \delta)n^q$ where $IC(K_q;V_1,...,V_q)=\prod\limits_{i<j}{d_{S}(V_i,V_j)}$. So to complete the proof it is enough for us to show that
$IC(B,W,\sigma)=IC(K_q;V_1,...,V_q)$. Indeed, we have
\begin{equation*}
IC(B,W,\sigma)=\prod\limits_{i<j}{d^{c(\sigma(i),\sigma(j))}(V_i,V_j)}=
\prod\limits_{i<j}{\frac{|E_{i,j}|}{|V_i||V_j|}}=\prod\limits_{i<j}{d_S(V_i,V_j)} =
IC(K_q;V_1,...,V_q)
\end{equation*}
\end{proof}

\begin{lemma}\label{aprox2}
For every $\delta>0$ and every $q$ there exists $\gamma=\gamma_{\ref{aprox2}}(\delta,q)$ such that the following holds: Suppose that $V_1,...,V_q$ are disjoint vertex sets of size $n$ each, and all pairs $(V_i,V_j)$ are $\gamma$-regular. Put $W=\{d^{\ell}(V_i,V_j) : 1 \leq i<j \leq q , ~1 \leq \ell \leq k\}$. Then for every $k$-colored graph $B$ on the vertices $[q]$, the number of copies of $B$ which have precisely one vertex in each of the sets $V_1,...,V_q$ is
\begin{equation*}
(IC(B,W) \pm \delta)n^q
\end{equation*}
\end{lemma}

\begin{proof}
Set $\gamma_{\ref{aprox2}}(\delta,q)=\gamma_{\ref{aprox1}}(\frac{\delta}{q!},q)$.
Let $V_1,...,V_q$ be as in the statement, and let $B$ be any $k$-colored graph. By Claim \ref{aprox1} for any permutation $\sigma: [q] \rightarrow [q]$, the number of copies of $B$ spanned by $v_1\in V_1,...,v_q\in V_q$ such that $v_i$ plays the role of $\sigma(i)$ is
$\left(IC(B,W,\sigma) \pm \frac{\delta}{q!}\right)n^q$. If we sum over all permutations
$\sigma: [q] \rightarrow [q]$, we get
$\sum\limits_{\sigma}{\left(IC(B,W,\sigma) \pm \frac{\delta}{q!}\right)}n^q$.
In this summation, we count every copy of $B$ exactly $Aut(B)$ times. Thus, by dividing by $Aut(B)$, we get that the number of copies of $B$ is
\begin{eqnarray*}
\frac{1}{Aut(B)}
\left( \sum\limits_{\sigma}{\left(IC(B,W,\sigma) \pm \frac{\delta}{q!}\right)}n^q \right) &=&
\left( \frac{1}{Aut(B)}\sum\limits_{\sigma}{IC(B,W,\sigma)} \pm \delta \right) n^q\\
&=& (IC(B,W) \pm \delta)n^q
\end{eqnarray*}
\end{proof}

\begin{lemma}\label{aprox3}
For every $\delta>0$ and $q$ there are $\gamma=\gamma_{\ref{aprox3}}(\delta,q)$ and $t=t_{\ref{aprox3}}(\delta,q)$ such that the following holds: Suppose that $R$ is a $k$-colored regularity instance of order at least $t$ and regularity parameter at most $\gamma$. Then for every $k$-colored graph $B$ on $q$ vertices, the number of copies of $B$ in any $n$-vertex $k$-colored graph satisfying $R$ is
\begin{equation*}
(IC(B,R) \pm \delta)\binom{n}{q}
\end{equation*}
\end{lemma}

\begin{proof}
Put
\begin{equation*}
t=t_{\ref{aprox3}}(\delta,q)=\left\lceil  \frac{4q^2}{\delta}\right\rceil
\end{equation*}
and
\begin{equation*}
\gamma=\gamma_{\ref{aprox3}}(\delta,q)=\min \left\{
\frac{\delta}{4q^2},
\gamma_{\ref{aprox2}}\left(\frac{\delta}{4},q\right)
\right\}\;.
\end{equation*}

Let $R$ be a $k$-colored regularity instance as in the statement, let $G$ be an $n$-vertex $k$-colored graph satisfying $R$ and let $B$ be any $k$-colored graph on $q$ vertices. Let $V_1,...,V_r$ be an equipartition of $V(G)$ satisfying $R$. Let $\mathcal{C}$ be the collection of all $q$-tuples that have at most one vertex in each of the sets $V_i$. By a union bound, the number of $q$-tuples that have more than one vertex in one of the sets $V_i$ is at most
\begin{equation*}
r\left(\frac{n}{r}\right)^2\binom{n-2}{q-2} \leq \frac{q^2}{r}\binom{n}{q} \leq \frac{q^2}{t}\binom{n}{q} \leq \frac{1}{4}\delta\binom{n}{q}\;.
\end{equation*}
So $|\mathcal{C}| \geq \left(1-\frac{\delta}{4}\right)\binom{n}{q}$. Therefore the lemma will follow from showing that the number of $q$-tuples belonging to $\mathcal{C}$ which span a copy of $B$ is
$\left(IC(B,R) \pm \frac{3}{4}\delta\right)|\mathcal{C}|$.

Given $A=\{x_1,...,x_q\} \subseteq \{1,...,r\}$ let $N(A)$ denote the number of $q$-tuples
$v_1 \in V_{x_1},...,v_q \in V_{x_q}$ which span a copy of $B$. We say that $A$ is {\bf good} if
all the pairs $(V_{x_i},V_{x_j})$ ($1 \leq i < j \leq q$) are $\gamma$-regular. Otherwise $A$ is called {\bf bad}. If $A$ is good we get from our choice of $\gamma$ via Lemma \ref{aprox2} that
\begin{equation*}
N(A)=\left(IC(B,W(A)) \pm \frac{1}{4}\delta\right)\left(\frac{n}{r}\right)^q
\end{equation*}
where we set
\begin{equation*}
W(A)=\{d^{\ell}(V_i,V_j) : i,j \in A, ~1 \leq \ell \leq k\}.
\end{equation*}
We can thus estimate the number of $q$-tuples belonging to $\mathcal{C}$ which span a copy of $B$ by
\begin{equation*}
\sum\limits_{A \subseteq [r], |A|=q}{N(A)} =
\sum\limits_{A ~ is ~ good}
\left(\left(IC(B,W(A)) \pm \frac{1}{4}\delta\right)\left(\frac{n}{r}\right)^q \right)
+ \sum\limits_{A ~ is ~ bad}{N(A)} =
\end{equation*}
\begin{equation*}
\sum\limits_{A \subseteq [r], |A|=q}
\left(\left(IC(B,W(A)) \pm \frac{1}{4}\delta\right)\left(\frac{n}{r}\right)^q \right) +
\sum\limits_{A ~ is ~ bad}\left( N(A) - IC(B,W(A))\left(\frac{n}{r}\right)^q \right) =
\end{equation*}
\begin{equation}
\left(IC(B,R) \pm \frac{1}{4}\delta \right)\binom{r}{q}\left(\frac{n}{r}\right)^q +
\sum\limits_{A ~ is ~ bad}\left( N(A) - IC(B,W(A))\left(\frac{n}{r}\right)^q \right)
\end{equation}

Since $(V_i,V_j)$ is $\gamma$-regular for every $(i,j) \notin \bar{R}$ there are at most $\gamma\binom{r}{2}$ pairs $(V_i,V_j)$ which are not $\gamma$-regular. Therefore the number of bad sets $A \subseteq \{1,...,r\}$ is at most
\begin{equation*}
\gamma\binom{r}{2}\binom{r-2}{q-2} \leq \gamma q^2\binom{r}{q} \leq
\frac{1}{4}\delta\binom{r}{q}
\end{equation*}
Using the facts that $0 \leq IC(B,W(A)) \leq 1$ and
$0 \leq N(A) \leq \left(\frac{n}{r}\right)^q$ for every $A \subseteq [r]$, and the bound on the number of bad sets, we get that
\begin{equation*}
\left|
\sum\limits_{A ~ is ~ bad}\left( N(A) - IC(B,W(A))\left(\frac{n}{r}\right)^q \right)
\right|
\leq
\frac{1}{2}\delta \binom{r}{q} \left(\frac{n}{r}\right)^q
\end{equation*}
By plugging the above inequality in (4) we get that the number of $q$-tuples belonging to $\mathcal{C}$ which span a copy of $B$ is
$\left(IC(B,R) \pm \frac{3}{4}\delta \right)\binom{r}{q} \left(\frac{n}{r}\right)^q$. Observe that $|\mathcal{C}|=\binom{r}{q}\left(\frac{n}{r}\right)^q$. Therefore the number of those $q$-tuples is $\left(IC(B,R) \pm \frac{3}{4}\delta \right)|\mathcal{C}|$, as required.
\end{proof}

\begin{proof}[Proof (of Lemma \ref{aprox})]
Put $t=t_{\ref{aprox}}(\delta,q,k)=t_{\ref{aprox3}}(k^{-\binom{q}{2}}\delta,q)$,
$\gamma=\gamma_{\ref{aprox}}(\delta,q,k)=\gamma_{\ref{aprox3}}(k^{-\binom{q}{2}}\delta,q)$. Let $R$ be a regularity instance of order at least $t$ and regularity parameter at most $\gamma$ and let $G$ be a $k$-colored graph satisfying $R$. Let $B \in \mathcal{B}$. By our choice of
$\gamma$ and $t$ via Lemma \ref{aprox3}, the number of copies of $B$ in $G$ is $(IC(B,R) \pm k^{-\binom{q}{2}}\delta)\binom{n}{q}$.
Clearly $|\mathcal{B}| \leq k^{\binom{q}{2}}$, so the number of copies of graphs $B \in \mathcal{B}$ in $G$ is
\begin{equation*}
\sum\limits_{B \in \mathcal{B}}{\left( \left(IC(B,R) \pm k^{-\binom{q}{2}}\delta \right)\binom{n}{q} \right)} =
\left( \sum\limits_{B \in \mathcal{B}}{IC(B,R)} \pm \delta \right) \binom{n}{q}
\end{equation*}
\end{proof}

\subsection{Proof of Lemma \ref{differ}}\label{subsec:differ}

We will derive Lemma \ref{differ} from the following lemma.

\begin{lemma}\label{differ1}
For every $\delta$ and $q$ there is
$\lambda=\lambda_{\ref{differ1}}(\delta,q)$ such that the following holds:
Let $R$ and $M$ be $k$-colored regularity instances of order $r$, and densities
$ \{ \eta^{\ell}_{i,j} : 1 \leq i<j \leq r,~ 1 \leq \ell \leq k \} $ and
$ \{ \mu^{\ell}_{i,j} : 1 \leq i<j \leq r,~ 1 \leq \ell \leq k \} $ respectively. Let $B$ be a $k$-colored graph on the vertices $[q]$. Suppose that $\mu^{\ell}_{i,j}=\eta^{\ell}_{i,j} \pm \lambda$.
Then
$$
\left|IC(B,R)-IC(B,M)\right| \leq \delta
$$
\end{lemma}

\begin{proof}
Put
\begin{equation*}
\lambda=\lambda_{\ref{differ1}}=(\delta,q)= \frac{\delta}{2^{\binom{q}{2}}q!}.
\end{equation*}
Let $R,M$ be $k$-colored regularity instances as in the statement. Let $A =\{x_1,...,x_q\} \subseteq \{1,...,r\}$, and put
$$
W_{R}(A)=\{ \eta^{\ell}_{x_i,x_j} ~:~ 1 \leq i<j \leq q ~, 1 \leq \ell \leq k \}
$$
and
$$
W_{M}(A)=\{ \mu^{\ell}_{x_i,x_j} ~:~ 1 \leq i<j \leq q ~, 1 \leq \ell \leq k \}\;.
$$
Denote the color of $(i,j)$ in $B$ by $c(i,j)$. Let $\sigma: [q] \rightarrow [q]$ be a permutation. By Definition $\ref{icdef1}$ we have
\begin{equation*}
\left| IC(B,W_{R}(A),\sigma)-IC(B,W_{M}(A),\sigma) \right| =
\left| \prod\limits_{1\leq i<j \leq q}{\eta^{c(\sigma(i),\sigma(j))}_{x_i,x_j}} -
\prod\limits_{1\leq i<j \leq q}{\mu^{c(\sigma(i),\sigma(j))}_{x_i,x_j}} \right| =
\end{equation*}
\begin{equation*}
\left| \prod\limits_{1\leq i<j \leq q}{\eta^{c(\sigma(i),\sigma(j))}_{x_i,x_j}} -
\prod\limits_{1\leq i<j \leq q}{\left(\eta^{c(\sigma(i),\sigma(j))}_{x_i,x_j} \pm \lambda\right)} \right|
\end{equation*}
Opening the parentheses in the above product gives $2^{\binom{q}{2}}-1$ summands, all of which are multiples of $\pm \lambda$. Therefore
$\left| IC(B,W_{R}(A),\sigma)-IC(B,W_{M}(A),\sigma) \right| \leq \lambda 2^{\binom{q}{2}}$.
By Definition \ref{icdef2}, the triangle inequality and our choice of $\gamma$ we have
\begin{eqnarray*}
\left| IC(B,W_{R}(A))-IC(B,W_{M}(A)) \right| &=&
\left| \frac{1}{Aut(B)}\sum\limits_{\sigma}\left(IC(B,W_{R}(A),\sigma)-IC(B,W_{M}(A),\sigma)\right)
\right| \\&\leq& q!\lambda 2^{\binom{q}{2}} = \delta\;.
\end{eqnarray*}
By Definition \ref{icdef3} and the triangle inequality we have
\begin{equation*}
\left| IC(B,R)-IC(B,M) \right| = \left| \frac{1}{\binom{r}{q}}\sum\limits_{A \subseteq \{1,...,r\}, ~ |A|=q} \left(IC(B,W_{R}(A))-IC(B,W_{M}(A))\right) \right| \leq \delta\;,
\end{equation*}
as needed.
\end{proof}

\begin{proof}[Proof (of Lemma \ref{differ})]
Put $\lambda=
\lambda_{\ref{differ}}(\delta,q,k)=
\lambda_{\ref{differ1}}(k^{-\binom{q}{2}}\delta,q)$. By the choice of $\lambda$ via Lemma \ref{differ1} we get that $|IC(B,R)-IC(B,M)|\leq k^{-\binom{q}{2}}\delta$ for every $B \in \mathcal{B}$. Since $|\mathcal{B}| \leq k^{\binom{q}{2}}$ the triangle inequality thus gives (\ref{closeinstances}).
\end{proof}

\subsection{Proof of Lemma \ref{chop2}}\label{subsec:chop}

We will derive Lemma \ref{chop2} from the following lemma.

\begin{lemma}\label{chop1}
Let $U,V $ be disjoint vertex sets in a graph satisfying $|U|=|V|=n \geq n_{\ref{chop1}}(\gamma,k)$. Suppose the bipartite graph $(U,V)$ is $\gamma$-regular with $d(U,V)=\mu$. Let $\{\eta^{\ell} : 1 \leq \ell \leq k \} $ be nonnegative reals satisfying $\sum_{\ell=1}^{m}{\eta^{\ell}}=\mu$ and $\sum_{\ell=m+1}^{k}{\eta^{\ell}}=1-\mu$. Then there is a $(k,m)$-coloring of $(U,V)$ such that the resulting $k$-colored graph is $2\gamma$-regular and satisfies $d^{\ell}(U,V)=\eta^{\ell}$ for every $1 \leq \ell \leq k$.
\end{lemma}

For the proof of Lemma \ref{chop1} we need the following standard Chernoff-type inequality:

\begin{lemma}\label{chernoff}
Suppose $X_1,...,X_m$ are independent Boolean random variables and ${\mathbb P}(X_i=1)=p_i $. Let $E=\sum_{i=1}^{m}{p_i}$. Then ${\mathbb P}(|\sum_{i=1}^{m}{X_i}-E| \geq \delta m) \leq 2e^{-2\delta^{2}m}$.
\end{lemma}

\begin{proof}[Proof (of Lemma \ref{chop1})]
We will show that the edges between $U$ and $V$ can be colored with colors ${1,...,m}$ in a way that satisfies the requirements. The same argument can be used to color the non-edges with colors ${m+1,...,k}$. First assume that $\mu \leq \gamma$. If this is the case, just color any $\eta^{\ell}n^2$ of the edges between $U$ and $V$ with color $\ell$, for every $1 \leq \ell \leq m$. This way we made sure that $d^{\ell}(U,V)=\eta^{\ell}$. Let $U'\subseteq U , V' \subseteq V$ with $|U'|,|V'|\geq 2\gamma n$. Before the coloring we had $d(U',V')=\mu \pm \gamma$. Therefore after the coloring we have $0 \leq d^{\ell}(U',V') \leq \mu+\gamma$, so $|d^{\ell}(U',V')-\eta^{\ell}| \leq \mu+\gamma \leq 2\gamma$.

Assume from this point on that $\mu \geq \gamma$. For every edge $e \in E(U,V)$ roll a die with sides ${1,...,m}$ so that probability of side $\ell$ is $\frac{\eta^{\ell}}{\mu}$.
If the die falls on side $\ell$ then color $e$ with color $\ell$. Then the expected number of edges of color $\ell$ is $\eta^{\ell}n^{2}$.
By Lemma \ref{chernoff}, the probability that the number of edges colored by $\ell$ deviates from its expectation by more than $n^{\frac{3}{2}}$ is at most $2e^{-2 n/\mu} \leq
2e^{-2n}$. If $n \geq n_{\ref{chop1}}(\gamma,k)$ then this probability is less than $\frac{1}{4k}$. This means, by a union bound, that with probability at least $\frac{3}{4}$ the number of edges colored by $\ell$ is $\eta^{\ell} n^2 \pm n^{\frac{3}{2}}$ for every $1 \leq \ell \leq m$.

\begin{claim}
With probability at least $\frac{3}{4}$ all sets $U'\subseteq U, V' \subseteq V$ such that $|U'|,|V'| \geq 2 \gamma n $ satisfy $d^{\ell}(U',V')=\eta^{\ell} \pm  \frac{3\gamma}{2}$
\end{claim}

\begin{proof}
Let $U'\subseteq U, V' \subseteq V$ such that $|U'|,|V'| \geq 2 \gamma n $. The density of edges between $U'$ and $V'$ before the coloring is $\mu \pm \gamma$. Therefore the expected density of edges with color $\ell$ is $\frac{\eta^{\ell}}{\mu}(\mu \pm \gamma)=\eta^{\ell} \pm \gamma$. So it is enough to show that with probability at least $\frac{3}{4}$, there are no sets $U',V'$ and color $1 \leq \ell \leq m$ such that the density of edges of color $\ell$ between $U'$ and $V'$ deviates from its expectation by more than $\frac{\gamma}{2}$.

By Lemma \ref{chernoff}, the probability that the density of edges of color $\ell$ between $U'$ and $V'$ deviates from its expected value by more than $\frac{\gamma}{2}$ is at most
$2e^{-2(\frac{\gamma}{2})^2 |U'||V'|/d(U',V')} \leq 2e^{-\gamma^2 |U'||V'|/2}$. We assumed that $|U'|,|V'| \geq 2 \gamma n$, so this probability is at most $2e^{-\gamma^4 n^2/2}$. The number of choices of sets $U',V'$ as above is at most $2^{2n}$, and the number of colors is at most $k$, so by a union bound we get: The probability that there are sets $U',V'$ and a color $\ell$ , such that the density of edges between $U',V'$ with color $\ell$ deviates from its expectation by more than $\frac{\gamma}{2}$ is at most $k2^{2n}2e^{-\gamma^4 n^2/2}$. This expression is less than $\frac{1}{4}$ if $n$ is large enough, namely $n \geq n_{\ref{chop1}}(\gamma,k)$.
\end{proof}

Getting back to the proof of Lemma \ref{chop1} we see that so far we proved that with probability at least $\frac{1}{2}$ the following two conditions hold:
\begin{enumerate}
\item $d^{\ell}(U,V)=\eta^{\ell} \pm n^{-\frac{1}{2}}$ for every $1 \leq \ell \leq m$.
\item $d^{\ell}(U',V')=\eta^{\ell} \pm \frac{3\gamma}{2}$ for every $1 \leq \ell \leq m$ and every two sets $U' \subseteq U, V' \subseteq V$ of size at least $2\gamma n$.
\end{enumerate}

Now take a coloring that satisfies conditions 1 and 2. Let us write
$d^{\ell}(U,V)=\eta^{\ell} + \epsilon^{\ell}$ where
$\left| \epsilon^{\ell} \right| \leq n^{-\frac{1}{2}}$. Observe that
\begin{equation*}
\sum\limits_{\ell=1}^{m}{\left(\eta^{\ell} + \epsilon^{\ell}\right)}=
\sum\limits_{\ell=1}^{m}{d^{\ell}(U,V)} = \mu = \sum\limits_{\ell=1}^{m}{\eta^{\ell}}
\end{equation*}
Therefore $\sum\limits_{\ell=1}^{m}{\epsilon^{\ell}}=0$. We can change the colors of at most
$m n^{\frac{3}{2}} \leq k n^{\frac{3}{2}}$ edges to make sure that $d^{\ell}(U,V)$ is exactly $\eta^{\ell}$. For every $U'\subseteq U, V' \subseteq V$ of size at least $2\gamma n$ this final change changes $d^{\ell}(U',V')$ by at most $\frac{k}{(2\gamma)^{2}n^{\frac{1}{2}}}$ which is less than $\frac{\gamma}{2}$ if $n \geq n_{\ref{chop1}}(\gamma,k)$. So in the end we have $d^{\ell}(U',V')=\eta^{\ell} \pm 2\gamma$ as required.
\end{proof}

\begin{remark}
In fact, we could have proved the following stronger lemma: Let $U,V$ be disjoint sets in a graph, $|U|=|V|=n \geq n_{\ref{chop1}}(\gamma, \epsilon,k)$. Suppose $U,V$ is $\gamma$-regular and $d(A,B)=\mu$, and
let $\{\eta^{\ell} : 1 \leq \ell \leq k \} $ be nonnegative numbers satisfying $\sum_{l=1}^{m}{\eta^{\ell}}=\mu$ and $\sum_{l=m+1}^{k}{\eta^{\ell}}=1-\mu$ . Then there is a $(k,m)$-coloring of $U,V$ such that the resulting $k$-colored graph is $\gamma(1+\epsilon)$-regular and satisfies $d^{\ell}(U,V)=\eta^{\ell}$. The choice of $2 \gamma$ in Lemma \ref{chop1} and in Definition \ref{chopdef} is for convenience.
\end{remark}

\begin{proof}[Proof (of Lemma \ref{chop2})]
Put $n_{\ref{chop2}}(\gamma,t,k)=t \cdot n_{\ref{chop1}}(\gamma,k)$. Let $R$ be a graph regularity instance of order $r \leq t$ and regularity parameter $\gamma$. Let $R'$ be a $(k,m)$-chopping of $R'$ and let $\{ \eta^{\ell}_{i,j} : 1 \leq i<j \leq r,~1 \leq \ell \leq k \} $ be the densities of $R'$. Let $G$ be a graph with at least $n_{\ref{chop2}}(\gamma,t,k)$ vertices that satisfies $R$. Let $V_1,...,V_r$ be a $\gamma$-regular equipartition of $V(G)$ that corresponds to $R$. For every $1 \leq i \leq r$ we have $|V_i|\geq
\frac{n_{\ref{chop2}}(\gamma,t,k)}{r} \geq n_{\ref{chop1}}(\gamma,k)$. If $(i,j) \notin \bar{R}$ apply Lemma \ref{chop1} for $V_i,V_j$ and $\{ \eta^1_{i,j},...,\eta^k_{i,j}\}$. Color the rest of the edges and non-edges arbitrarily. The resulting $k$-colored graph satisfies $R'$.
\end{proof}

\end{document}